\documentstyle[psfig, 12pt]{amsart}
\title[a geometric application
of Janson's inequality]
{Uniformly distributed distances: A geometric application
of Janson's inequality}

\author{J\'anos Pach}

\thanks{Supported by NSF grant CCR-94-24398,
PSC-CUNY Research Award 663472 and OTKA-T-020914. This paper was
partly written while Pach was visiting MSRI
Berkeley, as part of the Convex Geometry program.}

\address{\hskip-\parindent
J\'anos Pach,
City College and Courant Institute, NYU,
 	251 Mercer Street, New York, NY 10012}
\email{pach@@cims6.nyu.edu}

\author{Joel Spencer}

\address{\hskip-\parindent
Joel Spencer,
Courant Institute, NYU,
        251 Mercer Street, New York, NY 10012}
\email{spencer@@cs.nyu.edu}

\newtheorem{theorem}{Theorem}[section]
\newtheorem{lemma}[theorem]{Lemma}

\begin{document}

\begin{abstract}
Let $d_1\leq d_2\leq\ldots\leq d_{n\choose 2}$ denote the
distances determined by $n$ points in the plane. It is shown
that $\min\sum_i (d_{i+1}-d_i)^2=O(n^{-6/7})$, where the
minimum is taken over all point sets with minimal distance
$d_1 \geq 1$. This bound is asymptotically tight.
\end{abstract}

\maketitle

\section{A problem of Erd\H{o}s on distance gaps} 
Consider $n$ points in the plane with
minimum distance at least one.  List the $m={n\choose 2}$
distances between them in increasing order: 
$d_1\leq d_2\leq\ldots\leq d_m$. The numbers $d_{i+1}-d_i$
will be called {\em distance gaps}. Erd\H os raised the 
following problem. Determine or estimate

\begin{equation}
f(n)= \min \sum_{1\leq i<{n\choose 2}}(d_{i+1}-d_i)^2, \label{eq1}
\end{equation}

\noindent where the minimum is taken over all $n$-element point sets
with minimum distance one. In particular, he asked if this 
sum can be made arbitrarily small. 

This choice of function to be minimized may at first appear to be
capricious.  Suppose, however, that we fix $d_1=1$ (which we may 
assume by a simple scaling argument) and the {\em diameter} $D=d_m$.  
Then if the numbers $d_i$ could be chosen arbitrarily, 
$\sum_{1\leq i<m}(d_{i+1}-d_i)^2$ would attain its minimum when
the $d_i$ are equally spaced. We have

\begin{equation}
\sum_{1\leq i<m}(d_{i+1}-d_i)^2\geq(m-1)^2(\frac{D-1}{m-1})^2=
\Omega(D^2/n^2). \label{eq2}
\end{equation}

The geometric constraints make it impossible to achieve perfectly
even spacing.

It was shown in [EMPS91] that the number of distances belonging to
the interval $[D-1,D]$ is at most $O(D^{3/2})$. So even if these
distances are evenly spaced, the gaps between them contribute at
least $\Omega(D^{-3/2})$ to the sum $\sum (d_{i+1}-d_i)^2$.
Combining this with (\ref{eq2}), we obtain

\begin{equation}
f(n)\geq\min_D\left(\Omega(D^2/n^2) + \Omega(D^{-3/2})\right)=
\Omega(n^{-6/7}). \label{eq3}
\end{equation}

\par Our main objective is to show that this bound is asymptotically
tight.

\medskip

\noindent {\bf Theorem 1.} {\em There exists a set of $n$ points in
the plane with minimum distance at least one such that the distances
$d_i$ determined by them satisfy}
$$\sum_{1\leq i<{n\choose 2}}(d_{i+1}-d_i)^2\leq cn^{-6/7},$$
{\em where $c$ is an absolute constant.}

\medskip

It follows from (\ref{eq3}) that the diameter $D$ of such a point set
must be approximately $n^{4/7}$ and that the interpoint distances 
must be fairly uniformly distributed in the interval $[1,D]$. Our
construction described in Section 4 is based on a {\em Poisson 
process}. For the analysis, we use the ``No Bonds Theorem'' 
(Section 3), which is a continuous version of Janson's Inequality
(Section 2). Sections 5--6 contain the details of the proof of
Theorem 1.

\section{Janson's inequality}

Let $X$ be a finite underlying set, and let $P$ be a random
subset of $X$, whose elements are chosen independently with
probability $\Pr[x\in P]=p_x$. Let $\{ S_i : i\in I\}$ be a
system of subsets of $X$, and let $A_i$ denote the event
that $S_i\subseteq P$. If $S_i\cap S_j=\emptyset$ then $A_i$
and $A_j$ are independent. Let
\[ \nu = \sum\Pr[A_i\wedge A_j],\]
where the sum is taken over all unordered pairs $i\not=j$
with $S_i\cap S_j\not=\emptyset$, and let
$$M=\prod_{i\in I}\Pr[\overline A_i]=\prod_{i\in I}(1-\Pr[A_i]).$$

\medskip

\noindent{\bf Janson's Inequality [J90].} {\em If $\Pr[A_i]
< \varepsilon$ for every $i\in I$, then} 
$$ M\leq \Pr[\wedge_{i\in I} \overline A_i]\leq 
Me^{\nu/(2-2\varepsilon)}.$$ 

\medskip

\par Let $G$ be a finite graph with vertex set $V(G)$ and
edge set $E(G)$. We apply Janson's Inequality in the
specific case when $X=V(G)$ and $\{ S_i : i\in I\}=E(G)$,
i.e., $|S_i|=2$ for every $i$. Then
$$ M=\prod_{\{x,y\}\in E(G)}(1-p_xp_y),$$
and $\nu$ is the expected number of ``vees'' (paths of
length two) in $G|_P$, the subgraph of $G$ induced by $P$.

\medskip

\noindent{\bf Corollary.} {\em Assume that for every edge
$\{x,y\}$ of a graph $G$, $p_xp_y\leq\varepsilon$. Then,}

$$ M\leq \Pr[ G|_P \mbox{\small{ is empty}}]\leq 
Me^{\nu/(2-2\varepsilon)}.$$

\medskip

Note that if $\max_{x\in V(G)}p_x$ is small, then $M$ can be
well approximated by 
\[ \prod_{\{x,y\}\in E(G)} e^{-p_xp_y} = e^{-\mu},\]
where $$\mu =\sum_{\{x,y\}\in E(G)} p_xp_y$$ is the expected
number of edges in $G|_P$. 

\section{Poisson processes}

Now we extend the Corollary to the continuous cased to be
used in the sequel.

Let $X\subset$ I$\!$R$^2$ be a bounded Jordan measurable set,
and let $\sim$ be a symmetric binary relation (graph) on $X$
such that $\{(x,y) : x\sim y\}$ is a Jordan measurable
subset of I$\!$R$^2\times$I$\!$R$^2$. If $x\sim y$ for some
$x,y\in X$, then they are said to form a {\em bond}.
Furthermore, let $\varphi$ be a countably additive finite
measure on the Borel subsets of $X$, defined by
$$\varphi(Y) = \int_Y \rho(x)dx,$$
where $\rho : X \rightarrow [0,+\infty)$, the {\em density
function} of $\varphi$, is Riemann integrable.

Let $P\subset X$ be a random multiset given by a {\em
Poisson process} associated with the measure $\varphi$. More
precisely, $P$ is a random variable whose values are almost
surely unordered $i$-tuples of $X$ (possibly with
repetition) for some non-negative integer $i$ such that

$$\Pr[|P|=i] = \frac{\varphi^i(X)}{i!}\cdot
e^{-\varphi(X)}\;\;\;(i=0,1,2,\ldots),$$

\noindent and for a fixed $i$, $P$ can be obtained by selecting $i$
points from $X$ independently with uniform distribution with
respect to $\varphi$ and disregarding their order. It is now
easy to check that for any Borel set $Y\subseteq X$,

\begin{equation}
\Pr[|P\cap Y|=i] = \frac{\varphi^i(Y)}{i!}\cdot e^{-\varphi(Y)}
\;\;\;(i=0,1,2,\ldots),\label{eq4}
\end{equation}
where $|P\cap Y|$ counts the number of points of $P$
belonging to $Y$ with multiplicities. In particular, the
expected value of $|P\cap Y|$ is equal to $\varphi(Y)$.
Moreover, if $Y_1$ and $Y_2$ are disjoint then $|P\cap Y_1|$
and $|P\cap Y_2|$ are independent random variables.

Let $B$ and $V$ denote the number of bonds and the number of
``vees'' formed by the elements of $P$, respectively, i.e.,

$$B = |\left\{ \{x,y\} : x,y\in P, x\sim y \right\}|,$$
$$V = |\left\{ (x,\{y,z\}) : x,y,z\in P, y\not=z, \mbox{ and
} x\sim y,
x\sim z \right\}|,$$

\noindent and set $\mu=E[B], \nu=E[V]$.

\medskip

\noindent {\bf Theorem 2. (No Bonds Theorem)} {\em Let
$P\subset X$ be a random multiset obtained by a Poisson
process associated with the measure $\varphi$, and let $B$
denote the number of bonds between elements of $P$.
\par Then, with the above notations and assumptions,}
$$e^{-\mu} \leq \Pr[B=0] \leq e^{-\mu +\nu}.$$

\medskip

\noindent {\bf Proof:} For a fixed $n$, place a mesh
$X_n=\{(i/n,j/n) : i,j\in$ Z$\!\!\!$Z$\}$ of sidelength $n$
in the plane, and let $P_n$ be a random multisubset of $X_n$
obtained by the Poisson process associated with the measure

$$\varphi_n(Y) = \sum_{\left({i \over n},{j \over n}\right)\in Y}
\varphi\left(\left[{i \over n},{i+1\over n}\right]\times\left[
{j\over n},{j+1\over n}\right]\right)$$

\noindent for any finite subset $Y\subseteq X_n$. If $Y$ consists of a
single element $x\in X_n$ then, by (\ref{eq4}),

$$\Pr[x\in P_n] =
\sum_{i=1}^{\infty}\frac{\varphi_n^i(x)}{i!}\cdot 
e^{-\varphi_n(x)}=1-e^{-\varphi_n(x)}.$$

\noindent Moreover, the events $x\in P_n$ are independent for all
$x\in X_n$.

Let $G_n$ denote the graph on the vertex set $X_n$, whose
two points $x,y\in X_n$ are joined by an edge if and only if
$x\sim y$. Furthermore, let $B_n$ denote the number of bonds
formed by elements of $P_n$, where each bond is counted only
{\em once}. It follows from the Corollary in the last
section that

\begin{equation}
M_n\leq \Pr[B_n=0]\leq M_ne^{\nu_n},\label{eq5}
\end{equation}
where
$$M_n=\!\prod_{\{x,y\}\in E(G_n)}\!\left(1-\Pr[x\in P_n]\Pr[y\in P_n]\right),$$

$$\nu_n=\!\!\!\sum_{\scriptsize\begin{array}{cc} (x,\{y,z\}): \\
xy,xz\in E(G_n), y\not=z \end{array}}
\!\!\!\Pr[x\in P_n]\Pr[y\in P_n]\Pr[z\in P_n].$$

Notice that, as $n$ tends to infinity,

$$\Pr[x\in P_n]=1-e^{-\varphi_n(x)}=1-\exp\left(-\int_{x+[0,1/n]^2}
\rho(x)dx\right)\rightarrow 0$$

\noindent uniformly for all $x$. Since $\rho$ is Riemann
integrable,

$$\lim_{n\rightarrow\infty} M_n =
\lim_{n\rightarrow\infty}\exp\left(-\left(1+o(1)\right)\cdot
\!\!\sum_{\{x,y\}\in E(G_n)}\!\varphi_n(x)\varphi_n(y)\right)$$
$$=\exp\left(-{1\over 2}\int_{x\in X}\int_{\scriptsize\begin{array}{cc}
y\in X \\
y\sim x \end{array}}\rho(x)\rho(y)dy dx\right)=e^{-\mu},$$

\noindent and

$$\lim_{n\rightarrow\infty} \nu_n =
\lim_{n\rightarrow\infty}\left(1+o(1)\right)\cdot
\!\!\!\sum_{\scriptsize\begin{array}{cc} (x,\{y,z\}): \\
xy,xz\in E(G_n), y\not=z \end{array}}\!\!\!
\varphi_n(x)\varphi_n(y)\varphi_n(z)$$
$$= {1\over 2}\int_{x\in X}\int_{\scriptsize\begin{array}{cc} y\in X \\
y\sim x \end{array}}\int_{\scriptsize\begin{array}{cc} z\in X \\
z\sim x \end{array}} \rho(x)\rho(y)\rho(z)dzdydx=\nu.$$

It remains to show that
$$\lim_{n\rightarrow\infty}\Pr[B_n=0]=
\lim_{n\rightarrow\infty}\Pr[B=0],$$
and then the theorem follows from (\ref{eq5}).

By definition, $\varphi_n(X_n)=\varphi(X)=\int_X\rho(x)dx$.
Hence, in view of (\ref{eq4}),

$$\Pr[B_n=0]=\sum_{i=0}^{\infty}\Pr\left[B_n=0 \mid
|P_n|=i\right]\Pr\left[|P_n|=i\right]$$
$$=\sum_{i=0}^{\infty}\Pr\left[B_n=0 \mid
|P_n|=i\right]\cdot\frac{\varphi^i(X)}{i!}\cdot
e^{-\varphi(X)},$$
$$\Pr[B=0]=\sum_{i=0}^{\infty}\Pr\left[B=0 \mid
|P|=i\right]\cdot\frac{\varphi^i(X)}{i!}\cdot e^{-\varphi(X)}.$$

It suffices to verify that for any fixed $i$,

\begin{equation}
\lim_{n\rightarrow\infty}\Pr\left[B_n=0 \mid |P_n|=i\right]=
\Pr\left[B=0 \mid |P|=i\right]. \label{eq6}
\end{equation}

\noindent Restricting $P_n$ (and $P$) to the case
$|P_n|=i$ (and $|P|=i$), their points are distributed on
$X_n$ (on $X$) independently from each other and uniformly
with respect to $\varphi_n$ ($\varphi$, respectively). Thus,
$\Pr\left[B=0 \mid |P|=i\right]$ can be expressed as
$\int\rho(x_1)\ldots\rho(x_i)dx_1\ldots dx_i$ over a Jordan
measurable subset of $X\times\ldots\times
X\subseteq\; $I$\!$R$^{2i}$, and $\Pr\left[B_n=0 \mid
|P_n|=i\right]$ will approximate it with arbitrary
precision, as $n\rightarrow\infty$. This proves (\ref{eq6}),
and hence the theorem.
$\Box$

\section{Outline of the upper bound construction} 

The aim of this section is to sketch a probabilistic
construction for the proof of Theorem 1. The details will be
worked out in the next two sections. 

As we have indicated
in the last paragraph of Section 1, any $O(n)$-element point
set that satisfies the inequality in Theorem 1 must have
diameter $\Theta(n^{4/7})$. All of our points will be chosen
from the closed disk of radius $n^{4/7}$ around the origin 
$(0,0)\in$ I$\!\!$R$^2$.

\medskip

The construction consists of three parts (see Figure).

\bigskip
 
\centerline{\psfig{figure=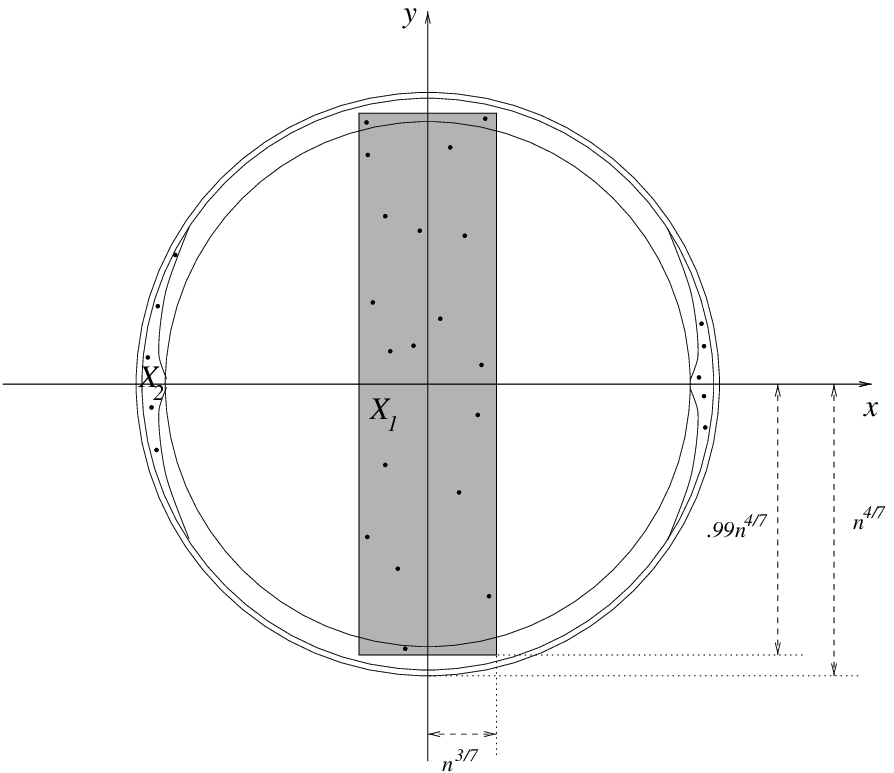}}
 
\centerline{\bf Figure}
 
\medskip

\noindent{\bf Step 1.} Define a point set $P_1$ by a Poisson
process on the rectangle
$$X_1=\{(x,y) : |x|\leq n^{3/7}, |y|\leq
0.99n^{4/7}\}$$
with constant density function
$\rho(x,y)\equiv\varepsilon$, where $\varepsilon$ is a small
positive number ($\varepsilon=10^{-3}$ will do). If two elements 
of $P_1$ are at distance less than one, delete both of them. 
(Subsection 5.1
and Section 6.)

\medskip

\noindent{\bf Step 2.} Define a point set $P_2$ by a Poisson
process with constant density $\varepsilon$ on the region
$X_2$ given by the polar coordinates
$$\{(r,\theta) : 0.9n^{4/7}<r<n^{4/7}-1, 
\min (|\theta|,|\theta -\pi|)<0.5(n^{4/7}-r)^{-1/4}\}.$$
As in the previous step, delete every pair of 
points whose distance is less than one. (Subsection 5.2.)

\medskip

\noindent{\bf Step 3.} We give an explicit construction
$P_3$ of points on the circle of radius $n^{4/7}$,
$$\{(x,y) : |x|^2+|y|^2 = n^{8/7}\},$$
with minimum distance at least one. (Subsection 5.3.) 

\bigskip

Note that the distance between any two points belonging to 
different $P_i$'s $(i=1,2,3)$ is larger than one.
We will show that the number of points in
$P=P_1\cup P_2\cup P_3$ is almost surely $O(n)$ and, with 
probability  at least 1/2, the sum of the squares of the distance
gaps determined by them satisfies
$\sum (d_{i+1}-d_i)^2=O(n^{-6/7}).$

First we need some terminology. Let us classify all
real numbers between 1 and $D=2n^{4/7}$ (the diameter of our
construction). A number $t$ is said to be
\[
\begin{array}{llrcl}
\mbox{{\it moderate  }} & \mbox{if  }& 1\leq&t&\leq 1.96n^{4/7}, \\
\mbox{{\it large  }} & \mbox{if  }& 1.96n^{4/7}<&t&\leq D-3,\\
\mbox{{\it extra large  }} & \mbox{if  }& D-3<&t&\leq D.
\end{array}
\]

For every  $j,\; 1\leq j<D=2n^{4/7}$ and for all $k=0,1,2,\ldots$, 
divide the interval $[j,j+1)$ into $2^k$ equal subintervals
of length $2^{-k}$. That is, let
\[ I_{jkl}=[j+(l-1)2^{-k},j+l2^{-k}),\;\;\;\;1\leq l\leq 2^k.
\]
In the sequel, these intervals will be called {\em
canonical}. Obviously, any interval $I\subset [1,D)$
with $|I|\leq 1$  has a canonical subinterval 
\[ I_{jkl}\subseteq I\;\; \mbox{ satisfying }\;\; |I_{jkl}|\geq
\frac{|I|}{4},\]
where $|I|$ denotes the {\em length} of $I$.

We say that an interval $I$ is {\em empty} if our
point set $P$ has no pair of points whose distance belongs
to $I$. Our goal is to show that with large
probability the sum of the squares of the distance gaps
determined by $P$ is small. We will see that our
construction has no distance gaps longer than one.
Since for any two consecutive distances, $d_i$ and
$d_{i+1}$, $(d_i, d_{i+1})$ is an empty interval, we have
\begin{equation}
\sum_{1\leq i<{n\choose 2}}(d_{i+1}-d_i)^2 \leq
4^2\cdot\!\!\!\!\sum_{\scriptsize\begin{array}{cc} I_{jkl} \mbox{
empty} \\ \mbox{canonical} \end{array}}|I_{jkl}|^2.\label{eq7}
\end{equation}
That is why in the rest of the paper we will concentrate
on establishing upper bounds for the expected value
\begin{equation}
E\:[\!\!\!\sum_{\scriptsize\begin{array}{cc} I_{jkl}
\mbox{ empty} \\ \mbox{canonical} \end{array}}|I_{jkl}|^2]
=\sum_{I_{jkl} \mbox{\scriptsize{ canonical}}}|I_{jkl}|^2 \Pr[I_{jkl}
\mbox{\small{ is empty}}]. \label{eq8}
\end{equation}

\section{Bounding the distance gaps}

\noindent{\bf 5.1. Moderate distances.} First we estimate 
the part of (\ref{eq8}), where
the sum is taken over all canonical intervals $I_{jkl}$ for
which $j$ is {\em moderate}, and the point
set is the set $P_1$ defined by a Poisson process in the
rectangle $X_1$, as described in Step 1 above. The expected
number of elements of $P_1$ is $\varepsilon\mbox{ Area}(X_1)
=3.96\varepsilon n=\Theta(n)$.     

Fix a canonical interval $I=I_{jkl}$ with
$1\leq j\leq 1.96 n^{4/7}$. We want to apply the
No Bond Theorem to the case when a pair of points $x,y\in
P_1$ forms a bond if  $|y-x|\in I$. Then the expected number
of bonds and the expected number of ``vees'' satisfy
\[ \mu={1\over 2}\int_{x\in X_1}\int_{\scriptsize\begin{array}{cc}
y\in X_1 \\ |y-x|\in I \end{array}}\varepsilon^2 dydx 
= \Theta(\varepsilon^2 n\min (j,n^{3/7})2^{-k}),\]
\[ \nu\leq \int_{x\in X_1}\int_{\scriptsize\begin{array}{cc}
y\in X_1 \\ |y-x|\in I \end{array}}\int_{\scriptsize
\begin{array}{cc}z\in X_1 \\ |z-x|\in
I\end{array}}\varepsilon^3 dzdydx=\Theta(\varepsilon^3
n\min (j^2,n^{6/7})2^{-2k}).\]
Thus, by the No Bond Theorem,
\[ \Pr [I \mbox{\small{ is empty}}]\leq e^{-\mu/2}\leq\exp\left(
-\Omega(\varepsilon^2 n\min (j,n^{3/7})2^{-k})\right),\]
whenever $|I|=2^{-k}\leq 1/\min (j,n^{3/7})$. In particular,
\begin{equation}
\Pr [I \mbox{\small{ is empty}}]\leq\exp\left(-\Omega(\varepsilon^2
n^{3/7})\right) \label{eq9}
\end{equation}
holds when $|I|=2^{-k}$ is roughly equal to $n^{-4/7}$,
and hence (\ref{eq9}) is valid for $k\leq (4/7)\log_2n$.

Therefore, the total contribution to (\ref{eq8}) of all
canonical intervals $I_{jkl}$ for which $j$ is moderate,

\[ \sum_{j \mbox{\scriptsize{ moderate}}}|I_{jkl}|^2
\Pr [I_{jkl} \mbox{\small{ is empty}}] \]
\[ =\sum_{j=1}^{1.96 n^{4/7}}\sum_{k=0}^{{4\over 7}\log n}
\sum_{l=1}^{2^k}2^{-2k}\Pr [I_{jkl}\mbox{\small{ is empty}}] \]
\[ +\sum_{j=1}^{1.96 n^{4/7}}\sum_{k>\frac{4}{7}\log n}
\sum_{l=1}^{2^k}2^{-2k}\Pr [I_{jkl}\mbox{\small{ is empty}}] \]
\[ \leq 1.96 n^{4/7}\left({4\over 7}\log n +1\right)
\exp\left(-\Omega(\varepsilon^2 n^{3/7})\right) \]
\[ + \sum_{j=1}^{1.96 n^{4/7}}\sum_{k\geq\frac{4}{7}\log n}
2^{-k}\exp\left(-\Omega(\varepsilon^2 n\min (j,n^{3/7})
2^{-k})\right) \]
\[ \leq \exp (-\Omega (n^{3/7})) + \sum_{j=1}^{1.96 n^{4/7}}
O\left(\frac{1}{\varepsilon^2 n\min (j,n^{3/7})}
\right)=O(n^{-6/7}). \]
Here we used the fact that the sum over $k$ can be majorated
by contsant times its largest term.

\medskip

\noindent{\bf 5.2. Large distances.} One can apply 
a similar argument to estimate the contribution of those
canonical intervals $I=I_{jkl}$ to (\ref{eq8}), for which 
$j$ is {\em large}. Now we have
to consider the Poisson process $P_2$ described in Step 2 
(see the previous section). Clearly, the
expected number of points of $P_2$ is equal to
\[ \varepsilon\mbox{Area}(X_2)=2\varepsilon\int_{0.9n^{4/7}}
^{n^{4/7}-1}0.5(n^{4/7}-r)^{-1/4}rdr=\Theta(\varepsilon
n).\] 

Fix a canonical interval $I=I_{jkl}$ with $1.96n^{4/7}<
j\leq D-3=2n^{4/7}-3$, and say that a pair of points forms a
{\em bond} if their distance belongs to $I$. Just like 
in the previous subsection, some routine calculation 
shows that the expected number of bonds is
\[ \mu=\Omega(\varepsilon^2 n^{6/7}(D-j)^{5/4}2^{-k}).\]
The No Bond Theorem now implies that
\[ \Pr [I \mbox{\small{ is empty}}]\leq e^{-\mu/2}\leq\exp\left(
-\Omega(\varepsilon^2 n^{6/7}(D-j)^{5/4}2^{-k})\right),\]
whenever $|I|=2^{-k}\leq n^{-4/7}$. This yields
\[ \Pr [I \mbox{\small{ is empty}}]\leq\exp\left(-\Omega(\varepsilon^2
n^{2/7})\right), \]
for every $k\leq (4/7)\log_2n$. 

Therefore, the total contribution to (\ref{eq8}) of all
canonical intervals $I_{jkl}$ for which $j$ is large,

\[ \sum_{j \mbox{\scriptsize{ large}}}|I_{jkl}|^2
\Pr [I_{jkl} \mbox{\small{ is empty}}] \]
\[\leq \sum_{j=1.96n^{4/7}}^{2n^{4/7}-3}
\sum_{k=0}^{{4\over 7}\log n}\sum_{l=1}^{2^k}2^{-2k}\Pr
[I_{jkl}\mbox{\small{ is empty}}] \]
\[+ \sum_{j=1.96n^{4/7}}^{2n^{4/7}-3}
\sum_{k\geq\frac{4}{7}\log n}\sum_{l=1}^{2^k}2^{-2k}\Pr
[I_{jkl}\mbox{\small{ is empty}}] \]
\[ \leq 0.04n^{4/7}\left({4\over 7}\log n +1\right)
\exp\left(-\Omega(\varepsilon^2 n^{2/7})\right) \]
\[+ \sum_{j=1.96n^{4/7}}^{2n^{4/7}-3}
\sum_{k\geq\frac{4}{7}\log
n}2^{-k}\exp\left(-\Omega(\varepsilon^2
n^{6/7}(D-j)^{5/4}2^{-k})\right) \]
\[ \leq \exp(-\Omega(n^{2/7})) + \sum_{j=1.96n^{4/7}}^{D-3}
O\left(\frac{1}{\varepsilon^2 n^{6/7}(D-j)^{5/4}}\right)
=O(n^{-6/7}).  \]

\noindent{\bf 5.3. Extra large distances.}
Here we present an explicit (i.e., non-probabilistic)
construction of a set $P_3$ of $O(n^{4/7})$ points on the
circle of radius $n^{4/7}$ centered at the origin $O$. Let
\[ P_3=\{ p_s : 0\leq s\leq n^{4/7}/2\} \cup
\{ q_t : 0\leq t\leq n^{4/7}/2\}, \]
where, using polar coordinates $(r,\theta)$,
\[ p_s=(n^{4/7}, 2sn^{-4/7}),\;\;\;
q_t=(n^{4/7}, \pi + 2t(n^{-4/7}+4n^{-8/7})).\]
For any $t\geq s$, the clockwise angle $p_sOq_t$ is
\[ \pi - 2(t-s)n^{-4/7} - 8tn^{-8/7}.\]
Thus, the angles $p_sOq_t\;\;(0\leq s\leq t)$ are
fairly densely distributed in the interval $(\pi - 
1/2,\pi)$; every closed subinterval of length $8n^{-8/7}$
contains at least one of them. Consider the distances
corresponding to these angles. It is easy to see that they
divide $[2n^{4/7}-3,\;2n^{4/7}]$ into
subintervals of length at most $15n^{-6/7}$. Consequently,
the sum of the squares of the distance gaps in this interval
is at most $45n^{-6/7}$, and
\begin{equation}
\sum_{\scriptsize\begin{array}{cc} I_{jkl} \mbox{
empty} \\ j \mbox{ extra large} \end{array}}|I_{jkl}|^2
= O(n^{-6/7}).  \label{eq10}
\end{equation}

\section{Excluding point pairs at distance less than one}

So far we have focused on bounding (\ref{eq8}), the expected
value of the sum of the squares of the distance gaps
determined by our random construction $P_1\cup P_2\cup P_3$.
The deterministic part of the construction, $P_3$, satisfies
the requirement that it has no two points at distance less than
one, and it is also true that the distance between any two
points belonging to different $P_i$'s is at least one.
However, some of the distances induced by $P_1$ and $P_2$ 
may be shorter than one.

In this section, we will deal with this problem. Let $P_i^*$
denote the point set obtained from $P_i$ by deleting every
point $x\in P_i$ for which there is another point $y\in P_i$
with $|y-x|<1\;\;(i=1,2)$. In fact, instead of (\ref{eq8}),
we need an upper bound on
\[ \sum_{I_{jkl} \mbox{\scriptsize{ canonical}}} |I_{jkl}|^2
\Pr{^*}[I_{jkl} \mbox{\small{ is empty}}],\] 
where $\Pr^*$ denotes the probability that a given
condition is satisfied for $P_1^*\cup P_2^*\cup P_3$, i.e.,
after the deletions have been carried out.

Let us concentrate on {\em moderate} distances, i.e., 
on the set $P_1^*$. Fix again a canonical interval $I=I_{jkl}$ 
with $1\leq j\leq 1.96 n^{4/7}$. A pair of points $\{x,y\}
\subseteq P_1$ is said  to form a bond if  $|y-x|\in I$. A
bond is called {\em bad} if there exists $z\in P_1$ whose
distance from $x$ or from $y$ is less than one. In this
case, $\{x,y\}\not\subseteq P_1^*$. Two bonds $\{x_1,y_1\}$
and $\{x_2,y_2\}$ are {\em separated} if their distance
is at least 2.  Let $\mu$ and $\nu$ denote the same as in 5.1.

\begin{lemma}\label{L1}
If $k\geq(4/7)\log_2n$, then the probability that there are 
at least $\mu/10$ pairwise
separated bonds is larger than $1-e^{-\mu/8}$.
\end{lemma}

\noindent{\bf Proof:}
We prove the stronger statement that the probability that 
there is a {\em maximal} family ${\cal F}$ of pairwise 
separated bonds with $|{\cal F}|<\mu/10$ is smaller than
$e^{-\mu/8}$.

Indeed, this probability is at most
\[ \sum_{i<{\mu\over 10}}\sum_{|{\cal F}|=i}\Pr[{\cal F}
\mbox{\small{ is maximal family of separated bonds}}]. \]
Since the expected number of $i$-tuples of bonds is
$\mu^i/i!$, the above sum cannot exceed
\begin{equation}
\sum_{i<{\mu\over 10}}\frac{\mu^i}{i!}\cdot
\max_{|{\cal F}|<{\mu\over 10}}\Pr [{\cal F}\mbox
{\small{ is maximal}} \; |\; {\cal F}\mbox{\small
{ consists of bonds}}].\label{eq11}
\end{equation} 
The conditional probabilities in the last expression can be
bounded by the No Bond Theorem. Delete from $X_1$ a unit
disk around  both points of each bond belonging to ${\cal
F}$, and restrict the Poisson process to the remaining set.
The probability that ${\cal F}$ is maximal is equal to the
probability that no bonds are formed under the restricted
process. Hence, this probability is at most $e^{-\mu'+\nu'}$,
where $\mu'<\mu$ and $\nu'<\nu$ are the expected number 
of bonds and ``vees'' in the restricted process,
respectively. If $k>(4/7)\log_2n$, then $\mu'>\mu/2>2\nu$.
This yields that the conditional probability in (\ref{eq11})
is at most $e^{-\mu'+\nu} < e^{-\mu/2+\nu} < e^{-\mu/4}$.
Therefore, (\ref{eq11}) can be bounded from above by
\[ \sum_{i<{\mu\over 10}}\frac{\mu^i}{i!}e^{-\mu/4}
<e^{-\mu/8},\]
as required.
$\Box$

\begin{lemma}\label{L2}
The probability that there are at least $\mu/10$ pairwise
separated bad bonds is at most $e^{-\mu/8}$.
\end{lemma}

\noindent{\bf Proof:} 
The expected number of {\em bad} bonds is at most
$2\varepsilon\pi\mu$, because for any bond, the probability 
that there is a point in the unit disk centered at
one of its points is $\varepsilon\pi$. If two point pairs are
separated, then the events that they form bad bonds are 
{\em independent}. Thus, the probability that there exist
$\mu/10$ bad bonds is at most 
\[ \frac{(2\varepsilon\pi\mu)^{\mu/10}}{(\mu/10)!}
<(100\varepsilon)^{\mu/10} < e^{-\mu/8}. \]
$\Box$

\medskip

It follows from Lemmas \ref{L1} and \ref{L2} that
if $j$ is {\em moderate} and $k\geq (4/7)\log_2n$, then the
probability that $P_1^{*}$ contains two points whose
distance belongs to $I_{jkl}$ satisfies
\[ \Pr{^*}[I_{jkl} \mbox{\small{ is empty}}]<2e^{-\mu/8}.\]
Thus, using the same estimates as in 5.1, we obtain
\begin{equation}
\sum_{j \mbox{\scriptsize{ moderate}}}\sum_{k\geq{4\over
7}\log n} \sum_{l=1}^{2^k}|I_{jkl}|^2
\Pr{^*}[I_{jkl} \mbox{\small{ is empty}}]=O(n^{-6/7}).
\label{eq12}
\end{equation}

On the other hand, the probability that there exist a
moderate $j$, $k\leq (4/7)\log_2n$, and $l$ such that the
distance of no two points of $P_1^{*}$ belongs to the
canonical interval $I_{jkl}$ is at most
\[ O\left(\sum_{j=1}^{1.96n^{4/7}}2^{-(4/7)\log n}
\Pr{^*}[I_{j\lceil (4/7)\log n\rceil l}\mbox{\small{ is
empty}}]\right)=\exp(-\Omega(n^{3/7})).\]

The above argument can also be applied to {\em large}
distances. Then,
\begin{equation}
\sum_{j \mbox{\scriptsize{ large}}}\sum_{k\geq{4\over
7}\log n} \sum_{l=1}^{2^k}|I_{jkl}|^2
\Pr{^*}[I_{jkl} \mbox{\small{ is empty}}]=O(n^{-6/7}),
\label{eq13}
\end{equation}
and the probability that there exist a {\em large $j$}, $k\leq
(4/7)\log_2n$, and $l$ such that the distance of no two 
points of $P_2^{*}$ belongs to the canonical interval 
$I_{jkl}$ is at most $\exp(-\Omega(n^{2/7}))$. 

\medskip

\noindent{\bf Summarizing:} With probability $1-o(1),
\;|P_1^*\cup P_2^*\cup P_3|\leq |P_1\cup P_2\cup P_3|=O(n).$
With probability at least
$1-\exp(-\Omega(n^{2/7}))$, for every canonical interval
$I_{jkl}$ for which $j$ is moderate or large and $k\leq
(4/7)\log_2n$, there is a pair of points in $P_1^*\cup
P_2^*\cup P_3$, whose distance belongs to $I_{jkl}$.
By (\ref{eq10}), (\ref{eq12}), and (\ref{eq13}), the
expected value of the sum of the squares of all other
empty canonical intervals satisfies
\[ E\:[\!\!\!\!\!\!\sum_{\scriptsize\begin{array}{ccc} I_{jkl}
\mbox{ empty} \\j \mbox{ moderate or large} \\ k\leq
(4/7)\log n \end{array}}\!\!\!\!\!|I_{jkl}|^2\;]+\sum_{\scriptsize
\begin{array}{cc} I_{jkl} \mbox{ empty} \\j 
\mbox{ extra large}\end{array}}\!\!\!|I_{jkl}|^2 = O(n^{-6/7}). \]
Hence, by Markov's inequality, with probability at least
1/2, the sum of the squares of {\em all} canonical intervals
will be $O(n^{-6/7})$. In view of (\ref{eq7}), this implies
that there exists a {\em specific} $O(n)$-element point set
$P_1^*\cup P_2^*\cup P_3$ with minimum distance one, for
which the sum of the squares of the distance gaps is
$O(n^{-6/7})$. This completes the proof of Theorem 1.

\end{document}